\documentclass{amsart}

\usepackage[utf8]{inputenc}
\usepackage{amsmath,amsfonts,amssymb,amsthm,enumerate}
\usepackage{pgf,tikz}
\usetikzlibrary{arrows}

\newtheorem{theorem}{\rm\bf Theorem}[section]

\newtheorem{corollary}[theorem]{\rm\bf Corollary}

\theoremstyle{definition}

\title{On Pasch's Axiom and Desargues' Theorem  in Busemann's work}
\author{Marc Troyanov}
\date{September 14, 2016}
\author{Marc Troyanov} 
\address{Section de Math{\'e}matiques,   
\'Ecole Polytechnique F{\'e}derale de Lausanne, station 8,
1015 Lausanne - Switzerland} 
\email{marc.troyanov@epfl.ch}

\begin{document}

\begin{abstract}
We discuss the role played by the techniques from the  ``foundations of geometry'', 
and in particular  by Desargues' Theorem, in the work of H. Busemann.  
This note, will be part of a  forthcoming edition of Busemann’s Collected papers.
\end{abstract}

\maketitle

\section{Introduction}

The subject of \emph{foundations of geometry}, or \emph{axiomatic geometry} as one would call it nowadays, has been very active during the last decades of the nineteenth  and the first half of the twentieth century. A major contributor is D. Hilbert and his work on the subject during the period 1891--1902 is recorded in   \cite{Hilbert2004}, see also \cite{Pambuccian2013}. His book \emph{Grundlagen der  Geometrie} had a tremendous influence on the subject
but a significant  number of other  mathematicians have contributed to the field including M. Pasch, G. Peano, E. H. Moore, F. Schur, O. Veblen,  A. N. Whitehead, R. Baer, E. Artin,  G. Birkhoff and F. Russel.
A pleasant introduction to the field as it stood in 1940 is the small book \cite{Robinson1940} by G. B. Robinson.

The aim of this field was to rigorously build the classical  Euclidean, hyperbolic and  elliptic geometries as well as the affine and projective geometries, on a simple set of basic geometric axioms. The methods drew  largely on the techniques and results of  projective geometry as developed in the first half of the nineteenth  century by  J.-V. Poncelet, M. Chasles,  A. F. Möbius and K. von Staudt.

In the 1920's the developments of topology and metric geometry, in particular by K. Menger,  brought a new viewpoint on the subject. The foundational question may be restated as the following problem: \emph{To characterize the Euclidean, hyperbolic and elliptic space among all abstract metric spaces}. 
This is the point of view adopted by  H. Busemann already in the 1930's and throughout his career.

The metric approach can no longer be considered as based on elementary axioms since it assumes set theory, and in particular the field of real numbers. 
One of the features of the metric viewpoint is that  Hilbert's  congruence,  order and completeness axioms reduce to elementary properties of the metric spaces.
Another feature is that the notion of \emph{angle} need not be discussed at all in a metric axiomatization of the classical geometries.

But a consequence is that new geometries appear.  Namely \emph{all metric spaces whose geodesics share the  properties of those of the classical geometries must be taken into account at some point in the theory}. The description of all these metric spaces is equivalent to Hilbert Problem IV and has been a central theme in Busemann's work. The most important examples  are  the Minkowski geometries (generalizing Euclidean geometry) and the Hilbert geometries (generalizing hyperbolic geometries).

Metric spaces whose geodesics are isomorphic to the system of lines in (a convex domain of) the affine or projective space are named \emph{Desarguesian spaces} by Busemann. They have  been introduced in 
his book \emph{Metric Methods in Finsler Spaces and in the Foundations of Geometry}  \cite{Busemann1942}  and the theory has been further developped
in the books \emph{The Geometry of Geodesics} and \emph{Recent synthetic differential geometry}  \cite{Busemann1955,Busemann1970}.
The theory of   Desarguesian spaces  is perhaps one of the most beautiful subjects  in Buseman's work.

Our goal in this short article is to explain Pasch's Axiom and Desargues' Theorem and their role in Busemann's theory of  Desarguesian spaces.

\section{Straight spaces and Pasch's axiom.}

The notion of straight space was introduced by Busemann in \cite[p. 38]{Busemann1955}. A convenient definition\footnote{In \cite[p. 72]{Busemann1942} Busemann has a slightly different definition  that also allows geodesics to be isometric to a circle.} is:

\smallskip

\textbf{Definition.} A \emph{straight space} is a  finitely compact metric space  $(X,d)$ such that any
pair of distinct points is joined by a unique geodesic and all  geodesics are isomorphic to the real line $\mathbb{R}$.
Recall that a metric space is \emph{finitely compact}, or \emph{proper} if every bounded subset is precompact.
 
Geodesics in  straight spaces are called  \emph{lines}. The \emph{segment} between two points $p$ and $q$ is 
the set of points $x$ such that $d(p,x)+d(x,q) = d(p,q)$. We denote this set by $[p,q]$.  

One of the first articles published by Busemann is his 1932 paper \cite{Busemann1932a}. It contains the following result:

\begin{theorem}\label{thHorospheres}
 For a straight space $(X,d)$ the following conditions are equivalent:
 \begin{enumerate}[(1)]
  \item The topological dimension of $X$ is $2$.
  \item Pasch's axiom is satisfied:  If the three points $a,b,c$ are not on a line, then any line containing a point $p\in [a,c]$
  also contains a point in $[a,b] \cup [b,c]$.
  \item $X$ is homeomorphic to $\mathbb{R}^2$.
\end{enumerate}
\end{theorem}
A straight space satisfying any of the above condition is called a \emph{straight plane} by Busemann.

We sketch the proof in order to illustrate Busemann's techniques. 
We first  prove the implication  (2) $\Rightarrow$ (3). Choose three non aligned points $a,b,c$ in $X$ and fix a point $p$ in the segment
$[a,c]$ different from $a$ and $b$. We then build a triangle $ABC$ in $\mathbb{R}^2$ isometric to $abc$ and we denote by $P$
the point in $[A,C]$ such that $\|P-A\| = d(a,p)$. 

We now associate to an arbitrary point $Q$ in $\mathbb{R}^2$, distinct from $P$, a point $q$ in $X$ as follows.
The line through $P$ and $Q$ meets $[A,B] \cup [B,C]$ at a point $S $. If $S \in [A,B]$, we denote
by $s$ the point on $ [a,b]$ such that $d(s,b) = \|S-B\| $ and if  $S \in [B,C]$ we choose
 $s$ to be the corresponding point in $[b,c]$.

Finally we choose  $q$ to be the unique point on the line through $[p,s]$ such that 
$$
 \frac{d(q,s)}{\|Q-S\|} =  \frac{d(q,p)}{\|Q-P\|} =  \frac{d(p,s)}{\|P-S\|}.
$$ 
We just constructed a map $Q \mapsto q = \phi(Q)$ from  $\mathbb{R}^2$ to $X$. This map is clearly injective and 
it is surjective because $X$ is assumed to satisfy Pasch's axiom.
It is not difficult to check that  the map $\phi$ is a homeomorphism (see \cite[Chap. 3 \S 2]{Busemann1942} or  \cite[Chap. 3 \S 10]{Busemann1955}).

To prove the implication  (1) $\Rightarrow$ (2), Busemann  starts with four points $a,b,c,p$ as above and constructs the corresponding map $\phi : \mathbb{R}^2 \to X$. If $\phi$ is not surjective then there exists a point $z \in X$ whose distance $\delta$ to the image 
$\phi(\mathbb{R}^2) \subset X$ is strictly positive. In this situation, Busemann proves that the sphere centered at $z$ with radius $\delta/2$ has dimension $\geq 2$, and
from this he concludes that the ball with same center and radius has dimension $\geq 3$. The argument is slightly involved, see 
 \cite{Busemann1932a} or \cite[p. 52--53]{Busemann1955} for the details.\footnote{Actually, in the book \cite{Busemann1955} a local version version of the statement is proved, but the argument is the same. This local version is a key step in Busemann's proof that a 2-dimensional $G$-space is a manifold.} The conclusion is that if $X$ is two-dimensional, then the map $\phi : \mathbb{R}^2 \to X$ is
associated to any triangle $abc$ and  any point $p\in [a,c]$ is surjective. This fact is clearly equivalent to Pasch's axiom.
 
 
\qed

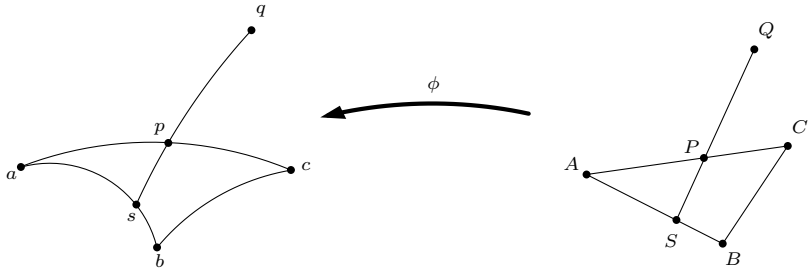
\begin{tikzpicture}[line cap=round,line join=round,>=triangle 45,x=0.5cm,y=0.5cm]
\clip(-11.2,-3) rectangle (22,7);
\fill[fill=black,fill] (-1.99,2.1) -- (-2.16,2.45) -- (-2.68,2.14)-- cycle;
\draw (4.4,0.62)-- (8.02,-1.22);
\draw (8.02,-1.22)-- (9.75,1.38);
\draw (9.75,1.38)-- (4.4,0.62);
\draw (6.79,-0.59)-- (8.86,3.95);
\draw [shift={(-9.87,-2.01)}] plot[domain=0.24:1.84,variable=\t]({1*2.93*cos(\t r)+0*2.93*sin(\t r)},{0*2.93*cos(\t r)+1*2.93*sin(\t r)});
\draw [shift={(-7.16,-8.14)}] plot[domain=1.18:1.94,variable=\t]({1*9.62*cos(\t r)+0*9.62*sin(\t r)},{0*9.62*cos(\t r)+1*9.62*sin(\t r)});
\draw [shift={(-2.47,-5.08)}] plot[domain=1.74:2.45,variable=\t]({1*5.9*cos(\t r)+0*5.9*sin(\t r)},{0*5.9*cos(\t r)+1*5.9*sin(\t r)});
\draw [shift={(8.47,-7.44)}] plot[domain=2.4:2.72,variable=\t]({1*17.61*cos(\t r)+0*17.61*sin(\t r)},{0*17.61*cos(\t r)+1*17.61*sin(\t r)});
\draw [shift={(0.3,-9.73)},line width=1.6pt]  plot[domain=1.36:1.8,variable=\t]({1*12.24*cos(\t r)+0*12.24*sin(\t r)},{0*12.24*cos(\t r)+1*12.24*sin(\t r)});
\begin{scriptsize}
\fill [color=black] (4.4,0.62) circle (1.5pt);
\draw[color=black] (4.0,0.9) node {$A$};
\fill [color=black] (8.02,-1.22) circle (1.5pt);
\draw[color=black] (8.3,-1.62) node {$B$};
\fill [color=black] (9.75,1.38) circle (1.5pt);
\draw[color=black] (10.08,1.89) node {$C$};
\fill [color=black] (7.52,1.06) circle (1.5pt);
\draw[color=black] (7.19,1.37) node {$P$};
\fill [color=black] (6.79,-0.59) circle (1.5pt);
\draw[color=black] (6.65,-1.15) node {$S$};
\fill [color=black] (8.86,3.95) circle (1.5pt);
\draw[color=black] (9.18,4.46) node {$Q$};
\fill [color=black] (-10.64,0.82) circle (1.5pt);
\draw[color=black] (-10.89,0.6) node {$a$};
\fill [color=black] (-7.57,-0.18) circle (1.5pt);
\draw[color=black] (-7.7,-0.5) node {$s$};
\fill [color=black] (-7.02,-1.32) circle (1.5pt);
\draw[color=black] (-6.95,-1.66) node {$b$};
\fill [color=black] (-6.72,1.46) circle (1.5pt);
\draw[color=black] (-6.95,1.83) node {$p$};
\fill [color=black] (-3.46,0.74) circle (1.5pt);
\draw[color=black] (-3.06,0.84) node {$c$};
\fill [color=black] (-4.51,4.46) circle (1.5pt);
\draw[color=black] (-4.26,4.97) node {$q$};
\draw[color=black] (0.33,3.02) node {$\phi$};
\end{scriptsize}
\end{tikzpicture}

\section{Desargues' Theorem}
\label{sec.desargues}

We now explain Desargues' Theorem. This theorem plays a basic role in Hilbert's foundations of Geometry
as well as in Busemann's work. 
The most convenient framework is that of elementary projective geometry. This is a subject where one considers only the undefined notions of points and lines\footnote{In this section and the next one, we follow the classical habit to note the point by capital latin letters and the line by smallcap letters} and the undefined notion of incidence. Intuitively, the point $P$ and the line $\ell$ are \emph{incident} if  $\ell$ passes through $P$.  Busemann and Kelly describe the subject in the following terms:  
\emph{Projective geometry is the totality of the facts which can be expressed solely in term of this incidence 
relation between points and lines}, see \cite[page 13]{BusemannKelly}.

This point of view, which essentially goes back  to Poncelet, 
is justified by the fact that a perspective projection from a plane $\pi$ to another plane $\pi'$
from a center $O$ maps points onto points, lines onto lines (with the exception of lines incident with $O$) and the incidence relations are preserved.
We refer to Chapter 3 of the book \cite{Hilbert1952} by Hilbert and Cohn-Vossen for an excellent concrete introduction to elementary projective geometry.

It may appear as a surprise that a rich geometry can be built on such an elementary principle.
The first non trivial result is probably   Desargues' Theorem:

\emph{Two triangles are such that the lines connecting the corresponding vertices are concurrent if and only if 
the intersection of the corresponding sides exist and are collinear.}

This means that if $ABC$ and $A'B'C'$ are two triangles such that the lines $AA'$, $BB'$ and $CC'$ meet at a point $O$, 
then the three points $I = AB\cap A'B'$, $J = AC\cap A'C'$ and $K = BC \cap B'C'$ exist and lie on a line $\ell$, and conversely.

\vspace{1cm}

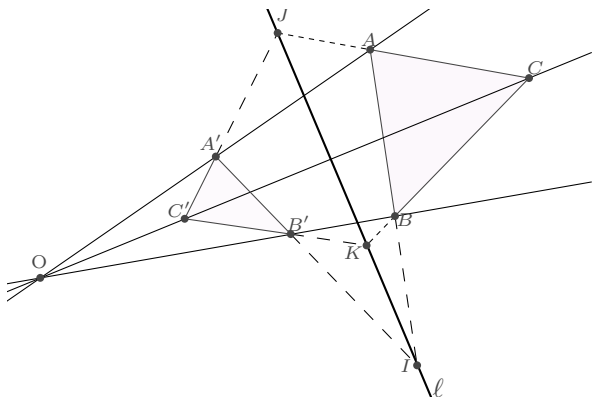
\definecolor{noir}{rgb}{0.27,0.27,0.27}
\definecolor{mauve}{rgb}{0.93,0.83,0.93}
\begin{tikzpicture}[x=0.6cm,y=0.6cm]
\clip(-1.2,-4.3) rectangle (15,4.3);
\fill[color=mauve,fill=mauve,fill opacity=0.15] (6.83,3.4) -- (7.38,-0.29) -- (10.35,2.77) -- cycle;
\fill[color=mauve,fill=mauve,fill opacity=0.1] (3.41,1.03) -- (5.07,-0.69) -- (2.72,-0.35) -- cycle;
\draw [color=noir] (6.83,3.4)-- (7.38,-0.29);
\draw [color=noir] (7.38,-0.29)-- (10.35,2.77);
\draw [color=noir] (10.35,2.77)-- (6.83,3.4);
\draw [domain=-1.3:11.74] plot(\x,{(-9.74--5.06*\x)/7.31});
\draw [domain=-1.3:11.74] plot(\x,{(-12.41--1.37*\x)/7.86});
\draw [domain=-1.3:11.74] plot(\x,{(-15.89--4.43*\x)/10.82});
\draw [color=noir] (3.41,1.03)-- (5.07,-0.69);
\draw [color=noir] (5.07,-0.69)-- (2.72,-0.35);
\draw [color=noir] (2.72,-0.35)-- (3.41,1.03);
\draw [dash pattern=on 2pt off 2pt] (4.78,3.77)-- (6.83,3.4);
\draw [dash pattern=on 5pt off 5pt] (3.41,1.03)-- (4.78,3.77);
\draw [dash pattern=on 5pt off 5pt] (5.07,-0.69)-- (7.87,-3.6);
\draw [dash pattern=on 5pt off 5pt] (7.38,-0.29)-- (7.87,-3.6);
\draw [dash pattern=on 5pt off 5pt] (5.07,-0.69)-- (6.75,-0.94);
\draw [dash pattern=on 2pt off 2pt] (7.38,-0.29)-- (6.75,-0.94);
\draw [line width=0.8pt,domain=-1.3:11.74] plot(\x,{(--29.92-4.7*\x)/1.98});
\begin{scriptsize}
\fill [color=noir] (6.83,3.4) circle (1.5pt);
\draw[color=noir] (6.8,3.63) node {$A$};
\fill [color=noir] (7.38,-0.29) circle (1.5pt);
\draw[color=noir] (7.6,-0.38) node {$B$};
\fill [color=noir] (10.35,2.77) circle (1.5pt);
\draw[color=noir] (10.49,3.01) node {$C$};
\fill [color=noir] (-0.48,-1.66) circle (1.5pt);
\draw[color=noir] (-0.5,-1.3) node {$\mathrm{O}$};
\fill [color=noir] (3.41,1.03) circle (1.5pt);
\draw[color=noir] (3.3,1.34) node {$A'$};
\fill [color=noir] (5.07,-0.69) circle (1.5pt);
\draw[color=noir] (5.24,-0.45) node {$B'$};
\fill [color=noir] (2.72,-0.35) circle (1.5pt);
\draw[color=noir] (2.6,-0.13) node {$C'$};
\fill [color=noir] (7.87,-3.6) circle (1.5pt);
\draw[color=noir] (7.62,-3.61) node {$I$};
\fill [color=noir] (4.78,3.77) circle (1.5pt);
\draw[color=noir] (4.89,4.18) node {$J$};
\fill [color=noir] (6.75,-0.94) circle (1.5pt);
\draw[color=noir] (6.46,-1.1) node {$K$};
\end{scriptsize}
\draw[color=noir] (8.34,-4.1) node {$\ell$};
\end{tikzpicture}
 \vspace{0.3cm}

One may call the point $O$ the \emph{center of perspectivity} of $A'B'C'$ relative to 
$ABC$ and $\ell$ the \emph{axis of perspectivity}. Desargues' Theorem may then be 
stated as follows: \emph{Two triangles admit a center of perspectivity if and only if they admit an
axis of perspectivity}.

This theorem bears a certain complexity. The Desargues configuration consists of 10 lines and 10 
points. Each point is incident to 3 lines and each line is incident to 3 points (the configuration is said to be
\emph{self-dual}).

Let us briefly sketch the  proof. One works in projective three-space and  first assumes that the two triangles $ABC$ and $A'B'C'$ lie in different planes $\pi$ and $\pi'$.

Suppose the two triangles $ABC$ and $A'B'C'$  are perspective from a point $O$, i.e. the point $O = AA' \cap BB' \cap CC'$  exists. 
This implies that the points $A'$ and $B'$ lie on the plane defined by $OAB$, therefore the lines $AB$ and $A'B'$ intersect at a point 
$I = AB\cap A'B'$. In particular, the point $I$ belongs to the intersection 
$\pi \cap \pi'$. This intersection is therefore not empty and it is not a plane since we assume $\pi \neq \pi'$. Therefore $\pi \cap \pi'$ is a line that we shall denote by $\ell$. The same argument shows that  the points $J = AC\cap A'C'$ and $K = BC \cap B'C'$ exist and lie on $\ell$.
We just proved that the points $I$, $J$ and $K$ exist and are collinear.

\smallskip

Suppose  conversely that the intersections $I=AB\cap A'B'$,  $J = AC\cap A'C'$ and $K = BC \cap B'C'$ exist; this implies in particular that the four points $A,A',B,B'$ are coplanar. Likewise the points $A,A',C,C'$ and $B,B',C,C'$ are coplanar. These three planes  have a unique common point $O$, therefore the three intersection lines $AA'$, $BB'$ and $CC'$ are incident with $O$. 

\smallskip

If the two triangles $ABC$ and $A'B'C'$ lie in a common plane $\pi$, then one chooses an auxiliary plane $\pi''$ such that $\pi''\cap \pi = \ell$ (the line $IJK$)
and one builds a (well chosen) triangle $A''B''C''$ in $\pi''$. Applying Desargues' Theorem twice, one may conclude. We refer to \cite{Robinson1940} or any book on synthetic projective geometry for the details.

\qed

\textbf{Remarks.}  Desargues' Theorem and its proof call for a number of remarks.  The first remark is that we have loosely 
assumed that any two lines in the same plane meet at a point. We have thus excluded the case of parallel lines. This is no
problem in the context of projective geometry, but in affine geometry, parallels do exist and Desargues' Theorem must be handled
with care. We will discuss below a  ``local version''  of Desargues'  Theorem due to Busemann that applies to Euclidean and 
Hyperbolic geometries as well.

The second remark is that the proof requires only a small number of basic properties of projective geometry (two points determine a unique line, 
etc.) and is  thus well adapted both to a naïve, concrete approach of projective geometry as well as to a formal, axiomatic treatment, see e.g. 
\cite{Robinson1940}.

The last remark is that the use of  the existence of at least a point outside a given plane in the proof is no accident. There are  geometries satisfying all the axioms of plane projective geometry in which Desargues' property fails. 
These planes are called ``non Desarguesian planes'' by Hilbert. We refer to \cite{Weibel}  for a recent survey of this topic.

\section{The harmonic conjugate}

The mid-point of a segment is a fundamental notion in Euclidean and  affine geometry. It is
however not an invariant notion under projections. In projective geometry, the 
useful corresponding concept is that of \emph{harmonic quadruple} and \emph{harmonic conjugate}.

Consider three distinct points $A,B$ and $P$ on a line $\ell$ and choose
 an arbitrary point $Z$ not
on that line. Choose a point $X$ on $AZ$ distinct from $A$ and $Z$ and let us denote by $Y$ and $W$ the points $Y=PX\cap BZ$
and $XB\cap YA$. The point $Q = \ell \cap WZ$ is called the \emph{harmonic conjugate} of $P$ with
respect to $A,B$, and the four aligned points $A,B,P,Q$ are said to be a \emph{harmonic quadruple}.
The configuration of the four point $X,Y,Z,W$ and the $6$ lines connecting them is called a 
\emph{complete quadrangle}.

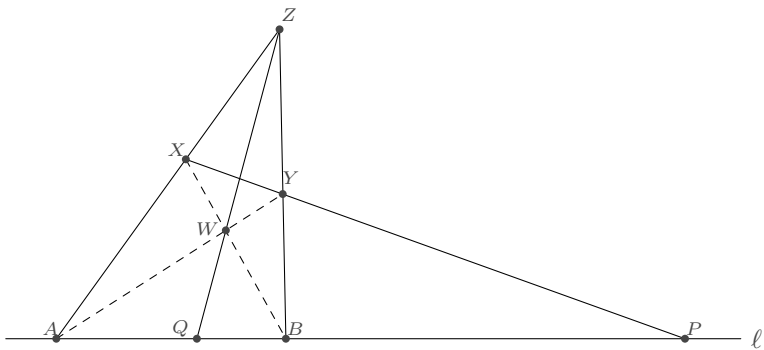
\begin{tikzpicture}[line cap=round,line join=round,>=triangle 45,x=0.7cm,y=0.7cm]
\clip(-1.96,-1.8) rectangle (14,6.96);
\draw [domain=-1.96:12] plot(\x,{(-0-0*\x)/1});
\draw (-1,0)-- (3.24,5.88);
\draw (3.36,0)-- (3.24,5.88);
\draw (1.46,3.41)-- (10.94,0);
\draw [dash pattern=on 3pt off 3pt] (-1,0)-- (3.3,2.75);
\draw [dash pattern=on 3pt off 3pt] (1.46,3.41)-- (3.36,0);
\draw (3.24,5.88)-- (1.67,0);
\draw[color=noir] (12.3,0) node {$\ell$};
\begin{scriptsize}
\fill [color=noir] (-1,0) circle (1.5pt);
\draw[color=noir] (-1.11,0.2) node {$A$};
\fill [color=noir] (3.36,0) circle (1.5pt);
\draw[color=noir] (3.54,0.2) node {$B$};
\fill [color=noir] (10.94,0) circle (1.5pt);
\draw[color=noir] (11.11,0.2) node {$P$};
\fill [color=noir] (3.24,5.88) circle (1.5pt);
\draw[color=noir] (3.42,6.15) node {$Z$};
\fill [color=noir] (1.46,3.41) circle (1.5pt);
\draw[color=noir] (1.29,3.6) node {$X$};
\fill [color=noir] (3.3,2.75) circle (1.5pt);
\draw[color=noir] (3.46,3.06) node {$Y$};
\fill [color=noir] (2.22,2.06) circle (1.5pt);
\draw[color=noir] (1.87,2.1) node {$W$};
\fill [color=noir] (1.67,0) circle (1.5pt);
\draw[color=noir] (1.36,0.2) node {$Q$};
\end{scriptsize}
\end{tikzpicture}

A fundamental consequence of Desargues' Theorem states that \emph{the harmonic conjugate $Q$ of
$P$ does not depend on the choice of the auxiliary point $Z$.}

\textbf{Proof.}   Choose a second auxiliary point $Z'$ not on the line $\ell$ and an arbitrary point $X'$
on $AZ'$  different from $A$ and $Z'$ and set $Y' = BZ\cap PX'$ and $W' = AY'\cap BX'$. We must prove that
the points $Q$, $Z'$ and $W'$ are aligned. This will follow from Desargues' Theorem applied thrice. 

Observe that by construction the points  $XY\cap X'Y'=P$, $XZ\cap X'Z=A'$ and $YZ\cap Y'Z'=B$ are aligned,
therefore, by Desargues' Theorem, the triangles $XYZ$ and $X'Y'Z'$ are  perspective from a point $O$, 
that is the point  $O= XX'\cap YY'\cap ZZ'$ exists.

By definition of $W$ and $W'$, we also have  $YW\cap Y'W'=A$, $XW\cap X'W'=B$
and $XY\cap X'Y'=P$,   aligned. Thus by Desargues' Theorem again the triangles $XYW$ and $X'Y'W'$ 
are perspective from a point, therefore the point $WW'\cap XX'$ coincides with $XX'\cap YY' = O$.

In particular the three lines $XX'$, $ZZ'$ and $WW'$ are incident with $O$. The triangles $XZW$ and
$X'Z'W'$ are thus  perspective and applying Desargues' Theorem a third time we conclude that $ZW \cap Z'W'$
is aligned with $YW\cap Y'W' = A$ and   $XW\cap X'W' = B$. 
This means that $ZW \cap Z'W'$ is on the line $\ell$ and, since $ZW \cap \ell = Q$, we have proved that 
$Z',W'$ and $Q$ are aligned.

\qed


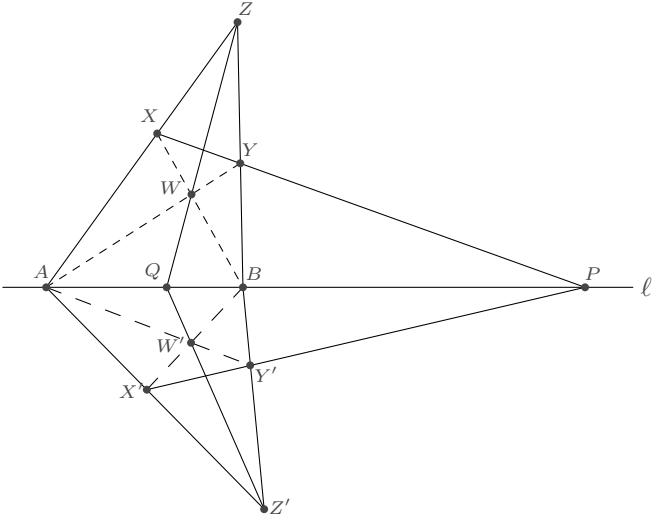
\begin{tikzpicture}[line cap=round,line join=round,>=triangle 45,x=0.6cm,y=0.6cm]
\clip(-1.96,-5.6) rectangle (14,6.5);
\draw [domain=-1.96:12] plot(\x,{(-0-0*\x)/1});
\draw (-1,0)-- (3.24,5.88);
\draw (3.36,0)-- (3.24,5.88);
\draw (1.46,3.41)-- (10.94,0);
\draw [dash pattern=on 3pt off 3pt] (-1,0)-- (3.3,2.75);
\draw [dash pattern=on 3pt off 3pt] (1.46,3.41)-- (3.36,0);
\draw (3.24,5.88)-- (1.67,0);
\draw (3.83,-4.92)-- (-1,0);
\draw (3.83,-4.92)-- (3.36,0);
\draw (3.83,-4.92)-- (1.67,0);
\draw (1.23,-2.27)-- (10.94,0);
\draw [dash pattern=on 6pt off 6pt] (-1,0)-- (3.52,-1.73);
\draw [dash pattern=on 6pt off 6pt] (1.23,-2.27)-- (3.36,0);
\draw[color=noir] (12.3,0) node {$\ell$};
\begin{scriptsize}
\fill [color=noir] (-1,0) circle (1.5pt);
\draw[color=noir] (-1.11,0.35) node {$A$};
\fill [color=noir] (3.36,0) circle (1.5pt);
\draw[color=noir] (3.6,0.3) node {$B$};
\fill [color=noir] (10.94,0) circle (1.5pt);
\draw[color=noir] (11.11,0.3) node {$P$};
\fill [color=noir] (3.24,5.88) circle (1.5pt);
\draw[color=noir] (3.42,6.18) node {$Z$};
\fill [color=noir] (1.46,3.41) circle (1.5pt);
\draw[color=noir] (1.29,3.8) node {$X$};
\fill [color=noir] (3.3,2.75) circle (1.5pt);
\draw[color=noir] (3.52,3.06) node {$Y$};
\fill [color=noir] (2.22,2.06) circle (1.5pt);
\draw[color=noir] (1.77,2.22) node {$W$};
\fill [color=noir] (1.67,0) circle (1.5pt);
\draw[color=noir] (1.36,0.35) node {$Q$};
\fill [color=noir] (3.83,-4.92) circle (1.5pt);
\draw[color=noir] (4.19,-4.85) node {$Z'$};
\fill [color=noir] (1.23,-2.27) circle (1.5pt);
\draw[color=noir] (0.9,-2.29) node {$X'$};
\fill [color=noir] (3.52,-1.73) circle (1.5pt);
\draw[color=noir] (3.88,-1.93) node {$Y'$};
\fill [color=noir] (2.21,-1.23) circle (1.5pt);
\draw[color=noir] (1.77,-1.27) node {$W'$};
\end{scriptsize}
\end{tikzpicture}

\section{Straight planes with Desargues' Property}

In this section we will sketch a proof of the following important Theorem of Busemann.

\begin{theorem} \label{ThDesargesPlanes}
 Let $(X,d)$ be a two-dimensional straight metric spaces satisfying the local Desargues property.
 Then   $X$ is projectively equivalent to a
 convex domain in  $\mathbb{R}^2$. That is, there is a map $\psi : X \to \mathbb{R}^2$ which 
 is injective and maps any geodesic in $X$ to an affine segment in  $\mathbb{R}^2$.
\end{theorem}

Note that the image $\psi(X) \subset \mathbb{R}^2$ is clearly a convex domain.

The ``local Desargues property'' referred to in this theorem is the statement that if two of the lines in a Desargues configuration meet at a point, then the third corresponding line also contains that point (remember that the classical Desargues configuration is a system of lines and points such that any point of the configuration belongs to three line). The precise formulation of the local Desargues property is given in 
\cite[pages 67--68]{Busemann1955}.

The proof of this theorem is quite involved. Here we only 
sketch the construction of the map $\psi$  and we refer to \cite[Section 13]{Busemann1955}
for the details.

\textbf{Step 1.}  Let us first consider three distinct points $a, b, p$ in  $X$ with $b \in [a,p]$
and choose
  a point $z$ not on the geodesic through $a$ and $b$. Choose
 also a point $x$ on the segment $[a,z]$ and construct the points
$y,w$ and $q$ as before. The existence of the points $x$ and $q$ follows from the two-dimensionality of $X$ and in particular from the 
validity of Pasch's axiom. 
By construction  $q$ is the harmonic conjugate of $p$ with respect to $a$ and $b$.

We   label the points $a,b,x,y,w,q$ as follows: $a(0,0)$, $b(1,0)$, $x(0,1)$, $y(1,1)$, $w(\frac{1}{2},\frac{1}{2})$
and $q(\frac{1}{2},0)$.

\smallskip 

\textbf{Step 2.} We  consider the points $s=pw\cap yb$, $t = zw \cap xy$ and $u=ax\cap wp$. The existence of these points is again
granted by Pasch's axiom. We label these new points as $s(1,\frac{1}{2})$, $t(\frac{1}{2},1)$ and $u(0,\frac{1}{2})$. 
The construction implies that $wpus$, $xypt$, $zwtq$ and $xauz$ are harmonic quadruples.

Let us call the procedure that defines the points $w,q,s,t,u$ the \emph{harmonic division} of the quadrilateral $abyx$.

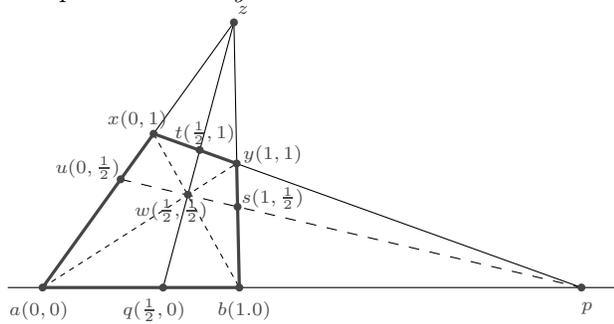
\begin{tikzpicture}[x=0.6cm,y=0.6cm]
\clip(-1.77,-1.46) rectangle (15,6.27);
\draw [domain=-1.77:11.78] plot(\x,{(-0-0*\x)/1});
\draw (-1,0)-- (3.24,5.88);
\draw (3.36,0)-- (3.24,5.88);
\draw (1.46,3.41)-- (10.94,0);
\draw [dash pattern=on 2pt off 2pt] (-1,0)-- (3.3,2.75);
\draw [dash pattern=on 2pt off 2pt] (1.46,3.41)-- (3.36,0);
\draw (3.24,5.88)-- (1.67,0);
\draw [dash pattern=on 4pt off 4pt] (0.73,2.4)-- (10.94,0);
\draw [line width=1.2pt,color=noir] (-1,0)-- (3.36,0);
\draw [line width=1.2pt,color=noir] (3.36,0)-- (3.3,2.75);
\draw [line width=1.2pt,color=noir] (3.3,2.75)-- (1.46,3.41);
\draw [line width=1.2pt,color=noir] (1.46,3.41)-- (-1,0);
\begin{scriptsize}
\fill [color=noir] (-1,0) circle (1.5pt);
\draw[color=noir] (-1.08,-0.5) node {$a(0,0)$};
\fill [color=noir] (3.36,0) circle (1.5pt);
\draw[color=noir] (3.47,-0.5) node {$b(1.0)$};
\fill [color=noir] (10.94,0) circle (1.5pt);
\draw[color=noir] (11.06,-0.45) node {$p$};
\fill [color=noir] (3.24,5.88) circle (1.5pt);
\draw[color=noir] (3.44,6.15) node {$z$};
\fill [color=noir] (1.46,3.41) circle (1.5pt);
\draw[color=noir] (1.1,3.72) node {$x(0,1)$};
\fill [color=noir] (3.3,2.75) circle (1.5pt);
\draw[color=noir] (4.1,2.95) node {$y(1,1)$};
\fill [color=noir] (2.22,2.06) circle (1.5pt);
\draw[color=noir] (1.85,1.7) node {$w(\frac{1}{2},\frac{1}{2})$};
\fill [color=noir] (1.67,0) circle (1.5pt);
\draw[color=noir] (1.47,-0.5) node {$q(\frac{1}{2},0)$};
\fill [color=noir] (0.73,2.4) circle (1.5pt);

\draw[color=noir] (0.,2.65) node {$u(0,\frac{1}{2})$};
\fill [color=noir] (3.32,1.79) circle (1.5pt);
\draw[color=noir] (4.1,1.98) node {$s(1,\frac{1}{2})$};
\fill [color=noir] (2.48,3.05) circle (1.5pt);
\draw[color=noir] (2.6,3.42) node {$t(\frac{1}{2},1)$};
\end{scriptsize}
\end{tikzpicture}

\textbf{Step 3.}  We iterate the procedure, that is,  we harmonically divide the four quadrilaterals $aqwu$, $qbsw$, $wsyt$ and 
$uwtx$. This produces 16 new points,   which we coherently label as $(\frac{m}{4}, \frac{n}{4})$ with $m,n \in \{0,1,2,4\}$. 

\begin{tikzpicture}[line cap=round,line join=round,>=triangle 45,x=0.9cm,y=0.9cm]
\clip(-1.6,-1) rectangle (15,6.33);
\draw [domain=-1.41:11.4] plot(\x,{(-0-0*\x)/1});
\draw (-1,0)-- (3.24,5.88);
\draw (3.36,0)-- (3.24,5.88);
\draw (1.46,3.41)-- (10.94,0);
\draw [dash pattern=on 2pt off 2pt] (-1,0)-- (3.3,2.75);
\draw [dash pattern=on 2pt off 2pt] (1.46,3.41)-- (3.36,0);
\draw (3.24,5.88)-- (1.67,0);
\draw [dash pattern=on 3pt off 3pt] (0.73,2.4)-- (10.94,0);
\draw [dash pattern=on 3pt off 3pt] (0.73,2.4)-- (2.48,3.05);
\draw [dash pattern=on 3pt off 3pt] (2.48,3.05)-- (3.32,1.79);
\draw [dash pattern=on 3pt off 3pt] (3.32,1.79)-- (1.67,0);
\draw [dash pattern=on 3pt off 3pt] (1.67,0)-- (0.73,2.4);
\draw [dash pattern=on 3pt off 3pt] (3.24,5.88)-- (0.5,0);
\draw [dash pattern=on 3pt off 3pt] (0.5,0)-- (2.6,0);
\draw [dash pattern=on 3pt off 3pt] (2.6,0)-- (3.24,5.88);
\draw [dash pattern=on 3pt off 3pt] (1.16,2.99)-- (10.94,0);
\draw [dash pattern=on 3pt off 3pt] (10.94,0)-- (0.09,1.51);
\draw [dash pattern=on 3pt off 3pt] (0.09,1.51)-- (1.16,2.99);
\begin{scriptsize}
\fill [color=noir] (-1,0) circle (1.2pt);
\draw[color=noir] (-1.07,-0.3) node {$(0,0)$};
\fill [color=noir] (3.36,0) circle (1.2pt);
\draw[color=noir] (3.44,-0.3) node {$(1,0)$};
\fill [color=noir] (10.94,0) circle (1.2pt);
\draw[color=noir] (11.03,-0.3) node {$P$};
\fill [color=noir] (3.24,5.88) circle (1.2pt);
\draw[color=noir] (3.41,6.03) node {$Z$};
\fill [color=noir] (1.46,3.41) circle (1.2pt);
\draw[color=noir] (1.37,3.66) node {$(0,1)$};
\fill [color=noir] (3.3,2.75) circle (1.2pt);
\draw[color=noir] (3.65,2.92) node {$(1,1)$};
\fill [color=noir] (2.22,2.06) circle (1.2pt);
\draw[color=noir] (1.75,2.15) node {$(\frac{1}{2},\frac{1}{2})$};
\fill [color=noir] (1.67,0) circle (1.2pt);
\draw[color=noir] (1.6,-0.3) node {$(\frac{1}{2},0)$};
\fill [color=noir] (0.73,2.4) circle (1.2pt);
\draw[color=noir] (0.3,2.45) node {$(0,\frac{1}{2})$};
\fill [color=noir] (3.32,1.79) circle (1.2pt);
\fill [color=noir] (2.48,3.05) circle (1.2pt);
\draw[color=noir] (2.4,3.3) node {$(\frac{1}{2},1)$};
\fill [color=noir] (1.8,2.8) circle (1.2pt);
\fill [color=noir] (2.72,1.14) circle (1.2pt);
\fill [color=noir] (0.5,0) circle (1.2pt);
\draw[color=noir] (0.5,-0.3) node {$(\frac{1}{4},0)$};
\fill [color=noir] (2.6,0) circle (1.2pt);
\draw[color=noir] (2.6,-0.3) node {$(\frac{3}{4},0)$};
\fill [color=noir] (0.09,1.51) circle (1.2pt);
\draw[color=noir] (-0.3,1.6) node {$(0,\frac{1}{4})$};
\fill [color=noir] (1.16,2.99) circle (1.2pt);
\fill [color=noir] (1.14,1.37) circle (1.2pt);
\fill [color=noir] (2.01,1.25) circle (1.2pt);
\draw[color=noir] (1.1,1.1) node {$(\frac{1}{4},\frac{3}{4})$};
\fill [color=noir] (2.87,2.47) circle (1.2pt);
\draw[color=noir] (2.8,2.64) node {$(\frac{3}{4},\frac{3}{4})$};
\fill [color=noir] (2.81,1.92) circle (1.2pt);
\fill [color=noir] (3.31,2.34) circle (1.2pt);
\fill [color=noir] (3.34,1.06) circle (1.2pt);
\draw[color=noir] (3.48,0.83) node {$(1,\frac{1}{4})$};
\end{scriptsize}
\end{tikzpicture}
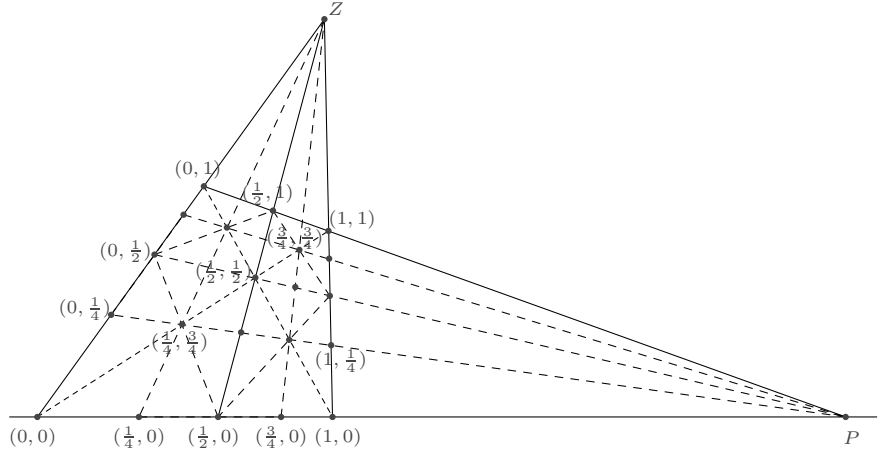
 \centerline{\footnotesize  Third step in the construction of the harmonic net with some of the labels.}

\bigskip

Repeating the procedure $k$ times, we produce new points in the convex hull of $abyx$ with label $(\alpha, \beta) = (\frac{m}{2^k}, \frac{n}{2^k})$
with    $0 \leq m,n \leq 2^k$.  
The totality of the  points thus constructed (as $k \to \infty$)  is called the \emph{harmonic net}  or the \emph{Möbius net} of the quadrilateral $abyx$. It can be shown to be a dense subset 
of the convex hull  $\mathcal{Q}$ of $abyx$. 

\smallskip

\textbf{Step 4.} We now build the map $\psi$. We first map any point $v$ 
of the harmonic net with dyadic label  $(\alpha, \beta) = (\frac{m}{2^k}, \frac{n}{2^k})$
to the point of $\psi(v) \in [0,1]^2 \subset \mathbb{R}^2$ with coordinates $(\alpha, \beta)$. This map can be shown to be continuous and can  therefore be extended to $\mathcal{Q}$. Furthermore the map $\psi : \mathcal{Q} \to [0,1]^2 $  can be proved to be a collineation, that is it maps geodesic segments in $\mathcal{Q}$ to affine segments in $[0,1]^2$.

\smallskip 
 
\textbf{Step 5.}  We finally extend the map as a global map $\psi : X \to \mathbb{R}^2$ as follows.
Given a point $v \in X\setminus \mathcal{Q}$, one considers two geodesics through $v$ containing interior points of $\mathcal{Q}$. These geodesics 
meet $\mathcal{Q}$ in segments $[s_1,s_2]$ and $[t_1,t_2]$. Let us denote by $t_1' = \psi(t_1)$,  $t_2' = \psi(t_2)$ and  $s_1' = \psi(s_1)$,  $s_2' = \psi(s_2)$ 
the images of these points in $[0,1]^2$.  We define $\psi(v) = v'$ to be the intersection of the lines $t_1't_2'$ and $s_1's_2'$. 
Using Desargues' Theorem again one may prove that  $\psi(v)$  does not depend on the chosen geodesics and that the globally defined map $\psi : X \to \mathbb{R}^2$ is a homeomorphism onto its image and maps geodesics of $X$ to affine segments of $\mathbb{R}^2$.

\qed

\section{Desarguesian spaces}

Theorem \ref{ThDesargesPlanes} has the following grand generalization:

\begin{theorem}\label{th.desspace}
Let $(X,d)$ be a $G$-space such that every pair of points can be joined by a unique geodesic. Assume either that $X$ has topological dimension $2$ and that the local Desargues property holds, or that $X$ has dimension $\geq 3$ and that any triple of points in $X$ belongs to a subspace that is itself a two-dimensional $G$-space for the induced metric. \\
Then either all geodesics of $X$ are topological circles and they have the same length and there is a homeomorphism from $X$ to $\mathbb{RP}^n$ that maps every geodesic in $X$ onto a projective line,
or $X$ is a straight space and there is a homeomorphism of $X$ to a convex domain in $\mathbb{R}^n$ that maps every geodesic in $X$ onto a line or an affine segment.
\end{theorem}

Recall that a $G$-space is a finitely compact metric space such that the geodesics can locally be extended
in a unique way. See \cite{Busemann1955} for the precise definition.

This theorem is proved in \cite[Section 14]{Busemann1955}. The hypothesis implies that the local Desargues property holds in $X$ for all dimension (the proof is a variant of the proof we gave in Section \ref{sec.desargues}).
This justifies the following definition:

\smallskip

\textbf{Definition.}  A metric space $(X,d)$ satisfying the hypothesis (and hence the conclusion) of  Theorem \ref{th.desspace}  is called  a \emph{Desarguesian space}.

\smallskip

One should note that the converse   is easy to prove: \emph{Any subset of $\mathbb{RP}^n$ that is either the whole of $\mathbb{RP}^n$ or a convex domain in $\mathbb{R}^n$ admits a complete metric whose geodesics are precisely the projective lines, respectively the affine segments. }

In the case of $\mathbb{RP}^n$ or $\mathbb{R}^n$, such a metric is given by the canonical spherical or Euclidean metric.
In the case of a convex domain  $\Omega$ in $\mathbb{R}^n$, then the following metric has the desired property:
$$
  d(x,y) = \|y-x\| + h(x,y)
$$
where $\|\; \|$ is a Euclidean metric and $h$ is the weak (degenerate) Hilbert metric in  $\Omega$ which is
defined as
\begin{equation}\label{weakHmetric}
  h(x,y) = \inf  \log\left(\frac{\|x-b\|}{\|y-b\|} \cdot \frac{\|y-a\|}{\|x-a\|}\right),
\end{equation}
where the infimum is taken over all pairs of points $a,b \in \Omega$ such that $y \in [x,b]$ and 
$x \in [a,y]$.\footnote{If the convex domain $\Omega$ is bounded, then (\ref{weakHmetric}) is a metric in the usual sense. In
general $h(x,y) = 0$ if and only if the entire line through $x$ and $y$ is contained in $\Omega$.}

The fourth Hilbert problem asks for a description of all Desarguesian spaces, and the  book \cite{Busemann1955} 
as well as later works by Busemann contains a rich supply  of properties of these spaces.

\textbf{Remark.}  Theorem \ref{th.desspace} has been extended to the (non metric) setting of  \emph{order geometry}
by J. P. Doignon and  J. Cantwell \& D. Kay  and   see \cite{Doignon} and \cite{Cantwell}.
See also the book \emph{Non-Euclidean Geometry} by Coxeter \cite{Coxeter}.
The book \cite{Coppel1998} by   W. A. Coppel calls it the \emph{Foundamental theorem of convex geometry}.

\section{Characterizations of the classical geometries}

We may now return to our original problem, that is \emph{to characterize the Euclidean, hyperbolic or elliptic space among all abstract metric spaces}. 
In the non compact case, one may then subdivide the problem in two subproblems:

 \begin{enumerate}[1.]
  \item To characterize  the Hilbert and Minkowski geometries among 
Desarguesian spaces.
  \item To characterize Euclidean geometry among Minkowski geometries and hyperbolic geometry among Hilbert geometries.
\end{enumerate}

The second problem has been thoroughly investigated by a number of mathematicians including Busemann himself, and a large supply of characterizations are available in the literature. 
We refer to  \cite{Amir} for the Euclidean case and to \cite{Guo2014} for the hyperbolic case.

The first problem has also been an important question in Busemann's work, see \cite{early} in this volume for a descriptions of Busemann's early work in this direction. We conclude by quoting the following result:

\begin{theorem}
 Among all non compact Desarguesian spaces,  the Hilbert and Minkowski geometries are  distinguished by the property that  an isometry of one geodesic onto another geodesic is a  projectivity.
\end{theorem}

This is Theorem 3 page 38 in  \cite{Busemann1970}.
One of the possible ways to define the notion of projectivity from a segment to another one is the condition that it is  a continuous map sending harmonic quadruples to harmonic quadruples. Equivalently this is a map obtained by composing  perspective projections.

Let us conclude by discussing some consequences of Busemann's work in Riemannian geometry. 
We first state the following  result, which is an immediate consequence of Beltrami's Theorem,
see  \cite[Theorem 15.3]{Busemann1955} or \cite{Matveev2006}:

\begin{theorem}
 A Riemannian manifold is a Desarguesian space if and only if it is isometric (after a possible rescaling of the metric)
 to the Euclidean space $\mathbb{R}^n$, the Hyperbolic space $\mathbb{H}^n$, or the projective space  $\mathbb{RP}^n$
 equipped with its standard metric.
\end{theorem}

This result provides a simple proof of the following classical fact, which is essentially due to F. Schur.\footnote{See \cite[\S 5--6]{Schur}, where Schur proves a slightly more general statement. See also \cite[chap. 1, Theorem 18]{SpivakIII} for a modern reference.}

\begin{corollary}
 Let $(M,g)$ be complete and simply connected Riemannian manifold without conjugate point. Assume that 
 $\dim(M) \geq 3$ and that any triple of points lie in a totally geodesic two-dimensional submanifold of $M$.
 Then $M$ has constant sectional curvature.
\end{corollary}

\textbf{Proof.}  The hypothesis that  $(M,g)$ is complete, simply connected and has no conjugate points 
means that it is a straight space by  Cartan-Hadamard's Theorem.  The hypothesis on triples of points implies
that $M$ is a Desarguesian space, we conclude by the previous Theorem.

\qed

\paragraph*{Acknowledgement.}  The author is thankful to  V. N. Berestovskii, V. Matveev, 
A. Papadopoulos and P. Petersen for their remarks on the manuscript.


\end{document}